\newcommand{\R}{\mathbb{R}}

\newcommand{\realn}{\mathbb{R}^n}
\newcommand{\realnum}{\mathbb{R}^{n+1}}

\documentclass{article}

\usepackage{amsmath}

\newtheorem{theorem}{\bf Theorem}
\newtheorem{lemma}{\bf Lemma}
\newtheorem{definition}{\bf Definition}
\newtheorem{example}{ \bf Example}
\newtheorem{remark}{\bf Remark}

\newcommand{\bt}[1]{\begin{theorem}
\label{t#1}\it }

\newcommand{\et}{\end{theorem} }

\newcommand{\bl}[1]{\begin{lemma}
\label{t#1}\it }

\newcommand{\el}{\end{lemma} }

\newcommand{\bd}[1]{\begin{definition}
\label{d#1}\it }

\newcommand{\ed}{\end{definition} }

\newcommand{\bp}[1]{\begin{example}
\label{p#1}\rm }

\newcommand{\ep}{\end{example}}

\newcommand{\br}[1]{\begin{remark}
\label{t#1}\it }

\newcommand{\er}{\end{remark} }

\usepackage{amssymb}

\begin{document}
\begin{center}

{\LARGE\bf Symmetries in some extremal problems between
two parallel hyperplanes
}


\bigskip

{\sc Monica Moulin Ribeiro Merkle}

\bigskip

\parbox{25cc}{
{\bf Abstract. }{\small Let $M$ be a compact hypersurface with boundary $\partial M=\partial D_1 \cup \partial D_2$, $\partial D_1 \subset \Pi _1$,
         $\partial D_2 \subset \Pi _2$, $\Pi_1$ and $\Pi _2$ two parallel hyperplanes in $\realnum$ ($n \geq 2$). Suppose that $M$ is contained in the slab
         determined by these hyperplanes and that the mean curvature $H$ of $M$ depends
         only on the distance $u$ to $\Pi _i$, $i=1,2$. We prove that these hypersurfaces are
         symmetric to a perpendicular orthogonal to $\Pi _i$, $i=1,2$, under  different conditions imposed on
         the boundary of hypersurfaces on the parallel planes: (i) when the angle of contact between $M$ and $\Pi _i$, $i=1,2$ is constant; (ii)
         when $\partial u / \partial \eta$ is a non-increasing function of the mean curvature of the boundary, $\partial \eta$  the inward normal; (iii)
         when $\partial u / \partial \eta$ has a linear dependency on the distance to a fixed point inside the body that hypersurface englobes; (iv)
         when $\partial D_i$ are symmetric to a perpendicular orthogonal to $\Pi _i$, $i=1,2$.
          } }

\end{center}

\renewcommand{\thefootnote}{}

\footnotetext{\it Submitted to Eletronic Research Announcements in Mathematical Sciences.}
\bigskip


\section{Introduction}

Let $\Pi _1$ and $\Pi _2$ be  parallel hyperplanes in $\realnum$ ($n \geq 2$). The slab determined by these hyperplanes is the set of
$\realnum$ with boundary  $\Pi _1 \cup \Pi _2$.
Let $\Omega$ be a connected, open and bounded subset of $\realnum$, contained in the slab,
such that $\partial \Omega = M \cup D_1 \cup D_2$, $D_1 = \overline{\Omega} \cap \Pi _1$, $D_2 = \overline{\Omega} \cap \Pi _2$ and $M$ is a compact hypersurface with boundary  $\partial M = \partial D_1 \cup \partial D_2$. Suppose that the mean curvature of $M$ depends
only on the distance to $\Pi _i$, $i=1,2$. We will prove that these hypersurfaces are symmetric to a perpendicular orthogonal to $\Pi _i$, $i=1,2$, under  different conditions imposed on the boundary of hypersurfaces on the parallel planes.

In the case $n=2$, $\Omega$ can be interpreted as the interior of a drop of liquid trapped between two parallel planes in the presence of gravity. In this physical problem, there are forces proportional to the area of the free surface, the surface tension and wetting energy proportional to the area on the plates wetted by the drop. We will assume that $D_{1}$
and $D_{2}$ are nonempty sets. The Euler-Lagrange equation for this problem implies that the mean curvature of the free surface depends on the distance
to the planes.

The key analytic tools used in the proofs of major theorems in this present work are presented in the section 2. They will alow us to compare the
hypersurface with itself, using a standard procedure, known as Alexandrov symmetrization and generalized by Serrin.

In section 3, we give a full proof of the following theorem:

\begin{theorem}
\label{tprincipal}
{\it
Let $\Pi _1$ and $\Pi _2$ be  parallel hyperplanes in $\realnum$.
Consider  an embedded compact connected $C^2$ hypersurface $M$ contained in the region of $\realnum$ between  $\Pi _1$ and $\Pi _2$, with $\partial M \cap \Pi _{i} \neq \emptyset$, $i=1,2$. Suppose that the mean curvature of $M$ depends
only on the distance to $\Pi _i$, $i=1,2$, and that the angle of contact between $M$ and $\Pi _i$, $i=1,2$, is constant
along $\partial M$ in each one of the two support hyperplanes (maybe with different constants). Then $M$ is rotationally
symmetric with respect to a perpendicular orthogonal to $\Pi _i$, $i=1,2$ .}
\et

This theorem has been  much used, but it seems that the proof has never been published in full,
being instead invariably left to a reader. Here, we give its proof in order to use it in later results
and make the presentation complete. We will use an  argument similar to the one that can be found at \cite{wente1}
for only one supporting hyperplane.

In section 4, we consider other symmetry results under different hypothesis on the boundary of hypersurface.


\section{The touching principle}

We begin with  some definitions and analytic lemmas that will be necessary for the proofs of results in the following sections.

Let $w(x)=w(x_1,...,x_n)$ be a differentiable function defined on some region of $\realn$. We will denote by $w_i=w_i (x_1,...,x_n)$ and
$w_{ij}=w_{ij} (x_1,...,x_n)$ the partial derivative of $w$ with respect to $x_i$ and second partial derivative of $w$ with respect to $i$ and $j$, respectively.
For $a_{ij}(x)$ and $b_i(x)$ being continuous functions in an open set $B\subset \realn$ satisfying $a_{ij}=a_{ji}$, $ 1 \leq i, j \leq n$, let $M(w)$
be a linear differential operator on $B$, defined by

\begin{equation}
\label{linear}
M(w) = \sum _{i,j=1}^n a_{i,j} (x) w_{ij} + \sum _{i=1}^n b_{i} (x) w_{i}.
\end{equation}

We say that $M(w)$ is elliptic on $B$ if

\begin{equation}
 \sum _{i,j=1}^n a_{i,j} (x) \xi_{i} \xi_{j} >0,
\end{equation}

for all $x \in B$ and for all $ \xi=(\xi _1, ...   , \xi _n) \neq (0,...,0)$.

We say that $M(w)$ is uniformly elliptic on $B$ with ellipticity constant $k>0$ if

\begin{equation}
 \sum _{i,j=1}^n a_{i,j} (x) \xi_{i} \xi_{j} \geq k|\xi|^2,
\end{equation}

for all $x \in B$ and for all $ \xi=(\xi _1, ...   , \xi _n) \neq (0,...,0)$.

The Hopf's maximum principle and Hopf's maximum principle on the boundary are well known results for studying behavior of linear elliptic operators and uniformly elliptic operators, respectively. However, the differential equation that the mean curvature satisfies is not linear, but quasi-linear. The following touching principles are going to play a substitute for our analysis, being a consequence  of the mentioned maximum principles.

\bl{intp} (Interior touching principle)
Let $L(w)$ be an operator
\begin{equation}
\label{quasilinear}
L(w)=M(w) + c(x)w,
\end{equation}
where $M(w)$ is an elliptic operator on $B$ as defined by \eqref{linear} and $c(x)$ a continuous function on $B$. If there is a function $w(x) \in C^2 (B)$ satisfying
\begin{enumerate}
  \item $L(w) \geq 0$ on $B$,
  \item $w(x) \leq 0$ on $B$,
  \item $w(x_0)=0$ for a $x_0 \in B$,
\end{enumerate}
then $w(x)\equiv 0$ on $B$.
\el

\bl{bp} (Boundary touching principle)
Let $B$ be a region in $\realn$ such that the boundary of $B$ in a neighborhood of  $x_0 \in \partial B $ is of class $C^1$. On  $\overline{B}$, we consider
an operator $L(w)$ of type \eqref{quasilinear}, where $M(w)$ is a uniformly elliptic operator  as defined by \eqref{linear} on $\overline{B}$ and $c(x)$  is
a continuous function on $\overline{B}$.
If there is a function $w(x) \in C^2 (B) \cap C^1 (\overline{B})$ satisfying
\begin{enumerate}
  \item $L(w) \geq 0$ on $B$,
  \item $w(x) \leq 0$ on $\overline{B}$,
  \item $\displaystyle{w(x_0)=0}$,
  \item $\displaystyle{\frac{\partial w}{ \partial \eta} (x_0)=0}$, where $\nu$ is the inward normal,
\end{enumerate}
then $w(x)\equiv 0$ on $B$.
\el

In the cases treated here,  the mean curvature $H$ of a hypersurface $M$ is a function of the height.
Suppose that  $M$ can be locally defined in a non-parametric way  by $x_{n+1} = u(x_1,...,x_n)$,
where $u$ is a smooth function in some bounded region $B$ of $\realn$. Hence, in $B$, $u$ satisfies the nonlinear
elliptic differential equation of second order

\begin{equation}
\label{media}
\text{div} \left(\frac{\nabla u}{\sqrt{1+ \lvert \nabla u \rvert ^2}} \right) = nH(x,u),
\end{equation}

that can be written as

\begin{equation}
\label{mediamaior}
\sum _{i,j=1}^n a_{i,j} (x,u,\nabla u) u_{ij}
=\frac{1}{W}  \sum _{i=1}^n  u_{ii} - \frac{1}{W^3} \sum _{i,j=1}^n  u_i u_j u_{ij}=nH(x,u),
\end{equation}

where $\displaystyle{W = \sqrt{1+ \lvert \nabla u \rvert ^2} }$. As a consequence of \eqref{mediamaior},

\begin{equation}
\label{mediaelitica}
\sum _{i,j=1}^n a_{i,j} (x,u,\nabla u) \xi _i \xi _j
=\frac{1}{W^3}\left[ (1+ \lvert \nabla u \rvert ^2) \lvert \xi \rvert ^2 - \sum _{i,j=1}^n  u_i u_j \xi_i \xi_j \right] \geq \frac{1}{W^3} \lvert \xi \rvert ^2.
\end{equation}

Consider now two hypersurfaces $M$ and $\overline{M}$ given nonparametrically by $x_{n+1} = u(x)$ and $x_{n+1} = \overline{u}(x)$ respectively, with
$x=(x_1,\dotsc,x_n) \in B$. Suppose that $u$ and $\overline{u}$ are $C^2$ on $B$ and satisfy the same mean curvature equation \eqref{media}. It is known that the difference $w(x)=u(x)-\overline{u}(x)$ is a solution to a linear partial differential equation $L(w)=0$, where $L$ is an operator as defined in \eqref{quasilinear}. To see this, we construct between $M$ and $\overline{M}$ a convex family of hypersurfaces $\{u^t\}$, given locally by
$u^t (x)=tu(x)+(1-t)\overline{u} (x)=\overline{u} (x)+tw(x)$, for $t\in [0,1]$ and $x \in B$. We introduce $H(t)=H(u^{t})$, where $H(t)$ is $C^{1}$ on $[0,1]$ and we calculate the integral of $H^{\prime} (t)$ in the interval $[0,1]$. By one side, we have
     \begin{equation*}
     n\int _{0} ^{1} H^{\prime} (t)\;dt=n [H(1)-H(0)]=n[H(u)-H(\overline{u})].
     \end{equation*}
Applying the mean value theorem, there exists $\theta \in (0,1)$ such that $n[H(u)-H(\overline{u})] = n\frac{\partial H}{\partial u} (\overline{u}+\theta w) w$.
By the other side, considering the same notation as before, since $u_{i} ^{t}=\overline{u}_{i}+tw_{i}$ and $u_{ij} ^{t}=\overline{u}_{ij}+tw_{ij}$ for
$i,j=1,\dotsc,n$, we have

    \begin{equation*}
    n\int_{0}^{1}\;H^{\prime} (t)\;dt=
    \int_{0}^{1}\;\frac{d}{dt}\left[\sum _{i=1}^n  \frac{u^t_{ii}}{(1+\lvert \nabla u^t \rvert ^2)^{1/2}} -
                                \sum _{i,j=1}^n \frac{ u^t_i u^t_j u^t_{ij}}{(1+\lvert \nabla u^t \rvert ^2)^{3/2}} \right]\;dt
    \end{equation*}
    \begin{equation*}
    = \int_{0}^{1}\;\left[\sum _{i=1}^n  \frac{w_{ii}}{(1+\lvert \nabla u^t \rvert ^2)^{1/2}} -
                                \sum _{i,j=1}^n \frac{ w_i u^t_j u^t_{ij} + u^t_i w_j u^t_{ij} + u^t_i u^t_j w_{ij}}{(1+\lvert \nabla u^t \rvert ^2)^{3/2}} \right]\;dt
    \end{equation*}

     \begin{equation*}
    + \int_{0}^{1}\;\left[-\sum _{i=1}^n  \frac{u^t_{ii}(u^t_1 w_1+...+u^t_n w_n)}{(1+\lvert \nabla u^t \rvert ^2)^{3/2}} +
                                \sum _{i,j=1}^n \frac{3 u^t _i u^t_j u^t_{ij} (u^t_1 w_1+...+u^t_n w_n)}{(1+\lvert \nabla u^t \rvert ^2)^{5/2}} \right]\;dt\; .
    \end{equation*}

     This expression is nothing more than a nonlinear partial differential equation of second order in $w$, with integrals as coefficients. We can write it briefly as

     \begin{equation}
     \sum_{i,j=1}^n  A^{ij} (x,w(x),\nabla w(x)) w_{ij} + \sum _{i=1}^n B^{i} (x,w(x),\nabla w(x))w_{i} + C(x)w=0
     \end{equation}

where

     \begin{equation}
     A^{ij} (x,w(x), \nabla w(x))=\int _{0}^{1}\; a_{ij} (x,\overline{u} +tw, \nabla \overline{u} + \nabla w )\;dt\;.
     \end{equation}

     It turns out that the differential equation is of type \eqref{quasilinear}, with coefficients $A^{ij}$ independent of the set $u^t$, thus also independent of $w$, and $M(w)$ is elliptic on $B$. The fact that the coefficients $B_i$ are continuous and $B$ is a limited set of $\realn$ result on the uniform ellipticity of the equation on $\overline{B}$, with ellipticity constant

     \begin{equation*}
     k \geq \frac{1}{\text{max} ((1+\lvert \nabla u \rvert ^2)^{3/2}, (1+\lvert \nabla \overline{u} \rvert ^2)^{3/2})}\; .
     \end{equation*}

\br{intp}
Observe that we can redo all steps with a more general mean curvature, which depends on $u$ and $\nabla u$, and verifying that the prescribed mean equation will be an equation of type \eqref{quasilinear} with the same properties as obtained here.
\er
\section{Proof of Theorem 1}

%


     Let $L$ be a line in  $\Pi_1$ or $\Pi_2$ and $\Pi$ be a hyperplane perpendicular to $L$ and as a consequence also perpendicular to $\Pi_1$ and $\Pi_2$.
     As already mentioned, the hypersurfaces $M \cup D_1 \cup D_2$ compose the boundary of a body $\Omega$, which is compact. So, there exists a hyperplane parallel to $\Pi$, that is tangent to $M$ or touches $M$ on its boundary, such that $M$ is entirely on one side of
     $\Pi$. We denote such hyperplane by $\Pi (0)$.

     We begin to push the hyperplane $\Pi (0)$ along the line $L$, "entering" the body $\Omega$. After a displacement
     of $t$, we denote by $\Pi (t)$ the obtained hyperplane parallel to $\Pi (0)$. For each $t>0$ the hyperplane $\Pi (t)$ divides $\realnum$ in two
     components: $R (t)$, which contains $\Pi (0)$ and its complement $L (t)$.
     The portion of $M$ which is in the component  $R (t)$, will be denoted by $M (t)$.
     For small values of $t$, since $\Pi (0)$ is tangent to $M$, it is possible to reflect $M (t)$ with respect to $\Pi (t)$ obtaining
     $\overline{M} (t)$, which is entirely contained in the body $\Omega$ without any intersections. We observe that this process of
     reflection keeps invariant the two hyperplanes $\Pi _1$ and $\Pi _2$. It is also valid that the value of mean curvature of $\overline{M} (t)$,
     since it depends only on the distance to hyperplanes $\Pi _1$ and $\Pi _2$ agrees to the value of mean curvature of $M(t)$ at points of the
     same height in both hypersurfaces. And, since the body $\Omega$ does not need to be
     convex,  $\overline{M} (t)$ may have many connected components after reflection process.

     We continue this reflection process obtaining $\overline{M} (t)$ till we have a hyperplane $\Pi (t^*)$ where it happens the
     first contact between the reflected portion of $M$,   $\overline{M} (t^*)$, and the original portion of $M \subset L(t)$,
     with normals to $M$ and $\overline{M} (t^*)$ pointing to the same direction.
     This contact may be at an interior point $P$ or on the boundary at a point $P$ in one of the planes
     $\Pi _i$, $i=1,2$, i.e., one of the two conditions happen:

     \begin{enumerate}
       \item $\overline{M} (t^*)$ will be tangent to $M$ at an interior point $P$.
       \item $\overline{M} (t^*)$ will be tangent to $M$ at a point $P$ in $\partial M \cap \Pi _{i}$, $i=1$ or $i=2$.
     \end{enumerate}

     In the first case: $P$ is an interior point of $\overline{M} (t^*)$.

     We choose a coordinate system $(x_1,...,x_n,u)$ in $\realnum$, with the origin at $P$ such that the tangent space to $M$ at $P$, which coincides with the tangent space to $\overline{M} (t^*)$ at $P$ is $u=0$. We can also choose the direction of the $u$ axis to be directed to the interior of $\Omega$. This means that in a small neighborhood $B=\{x \in \realn ; |x| < r\}$ of the tangent space at $P$, $M$ and $\overline{M} (t^*)$ can be represented
     nonparametrically by two $C^2$ functions $u=u(x)=u(x_1,...,x_n)$ and $\overline{u}=\overline{u}(x)=\overline{u}(x_1,...,x_n)$, respectively. Both functions, $u(x)$ and $\overline{u} (x)$ satisfy equation of the mean curvature $H(x,u)$, coincide on the origin and $u(x) \leq \overline{u} (x)$ in this neighborhood. Then, $u(x)\equiv \overline{u} (x)$.

     In the second case: $P$ is a boundary point of $\overline{M} (t^*)$.

     By hypothesis, the angle of contact is constant along $\partial M$ and $\Pi _i$. So, $\overline{M} (t^*)$ will be tangent to $M$ at a point $P$. As before,
     we choose a coordinate system $(x_1,...,x_n,u)$ in $\realnum$, with the origin at $P$ such that the tangent space to $M$ at $P$, which coincides with the tangent space to $\overline{M} (t^*)$ at $P$ is $u=0$. We can also choose the direction of the $u$ axis to be directed to the interior of $\Omega$, the direction of $x_1$ to the interior of $M$, so that the tangent space to $\partial M$ at $P$ contained in $\Pi _1$ or $\Pi_2$ is given by $u=x_1=0$.

     We have $M$ a $C^2$ hypersurface with boundary, so in a neighborhood of $P$,   $M$ and $\overline{M} (t^*)$ can be represented
     nonparametrically by $C^2$ functions $u(x)$ and $\overline{u}(x)$, defined on domains $A_1 \subset \realn$ and $A_2\subset \realn$, respectively, that have the boundary of class $C^1$. Observe that, in the case where the angle of contact is $\pi / 2$, the neighborhood can be described as $A_1=A_2=B=\{x \in \realn ; |x| < r, x_1 > 0\}$. In the cases where the contact angle is different of $\pi /2$, we will have $A_1 \subset A_2$ or $A_2 \subset A_1 $. We consider $B=A_1 \cap A_2$. Both functions satisfy:
     \begin{itemize}
      \item At the origin $\displaystyle{u = \overline{u}} $ and $\displaystyle{\frac{\partial u}{\partial \eta}=\frac{\partial \overline{u}}{\partial \eta}}$, because $M$ and $\overline{M} (t^*)$ have the same tangent plane at $P$;
      \item On $\overline{B}$, $u(x)$ and $\overline{u} (x)$ satisfy equation of the mean curvature $H(x,u)$ and $u(x) \leq \overline{u} (x)$ because $\overline{M} (t)$ still lies in the interior of $\Omega$, for $t< t^*$.
     \end{itemize}

     Then, $u(x)\equiv \overline{u} (x)$.

     To complete the proof, we can observe that, picking up other directions for the line $L$, we find other hyperplanes of symmetry.
     The perpendicular orthogonal to $\Pi _i$, $i=1,2$, that we are searching is the intersection of all these hyperplanes of symmetry.

\section{Other boundary conditions}

The method of Alexandrov  can be used in other situations.

The first result shows how the symmetry of the boundary in a variational problem with fixed boundary on two
parallel hyperplanes implies symmetry of the hypersurface contained in the slab.

\bt{tbordofixo}
Let $\Pi _1$ and $\Pi _2$ be two parallel hyperplanes in $\realnum$. Let $\Omega$ be a connected, open and
bounded subset of $\realnum$ ($n \geq 2$), contained in the slab determined by $\Pi _1$ and $\Pi _2$, such that
 $\partial \Omega = M \cup D_1 \cup D_2$, $D_i = \overline{\Omega} \cap \Pi _i \neq \emptyset$, $i=1,2$ and $M$ is a compact embedded
 hypersurface with boundary  $\partial M = \partial D_1 \cup \partial D_2$.
Suppose $D_i$ are symmetric about a hyperplane $\alpha$ orthogonal to $\Pi_i$ and that $\partial D_i = \partial D_{i}^{+} \cup  \partial D_{i}^{-}$,
where $\partial D_{i}^{+}$ are graphs of  nonnegative $C^2$ functions $f_i$ defined on domains in $\Pi_i \cap \alpha$, which are positive in the
interior of the domains and zero on the border. The part $\partial D_{i}^{-}$ is the reflection of $\partial D_{i}^{+}$ in respect to $\alpha$.

Suppose that $M$ is a hypersurface of class $C^2$ in its interior and in the portion not touching $\alpha$.
Suppose that the mean curvature $H(x,u)$ of $M$ depends only on the distance to $\Pi _i$, $i=1,2$, being a $C^1$ function on $u$.
Then $\Omega$ is symmetric about $\Pi$.
\et

{\bf Proof. }
     Let $L$ be a line orthogonal to $\alpha$ and contained in $\Pi _1$ or $\Pi _2$. Let $\Pi$ be a hyperplane perpendicular to $L$ and as a consequence also perpendicular to $\Pi_1$ and $\Pi_2$ and parallel to $\alpha$.
     As in the previous theorem, the hypersurfaces $M \cup D_1 \cup D_2$ compose the boundary of a body $\Omega$, which is compact and
     exists a hyperplane parallel to $\Pi$ and $\alpha$, that is tangent to $M$ or touches $M$ on its boundary and $M$ is entirely on one side of
     $\Pi$. We denote such hyperplane by $\Pi (0)$. Using a pushing argument as in the proof of Theorem 1, and keeping the same notations, we find a
     first contact point  $P$ which can be in  the interior or on the boundary of $M$.

%

      What we have to show now is that $\overline{M} (t^*)$ is the hyperplane $\alpha$.
      Suppose that this does not happen and let us exam each case:

     In the first case, if $P$ is an interior point of $\overline{M} (t^*)$, we can repeat the same argument and obtain that $M$ is symmetric to respect to hyperplane $\overline{M} (t^*)$. In particular, $\partial D_i$ are symmetric to respect to $\overline{M} (t^*) \cap \Pi _i$, which contradicts the assumption made about $\partial D_i$, $i=1,2$.

     In the second case, if $P$ is a boundary point of $\overline{M} (t^*)$, since $\partial D_{1}^{-}$ is the reflection of $\partial D_{1}^{+}$ in respect to $\alpha$, we must have again $\overline{M} (t^*)$ and $\alpha$ the same hyperplane.

     This proves that the plane of symmetry has to be the same of symmetry of the boundary. \hfill $\square$

\medskip

As a consequence, if $\partial D_{i}$ are circles with centers in the same perpendicular orthogonal to $\Pi _ i$ (maybe with different radius), then $M$ is
rotationally symmetric with respect to this perpendicular.

\medskip

In the theorem of last section, the condition that the angle of contact between $M$ and $\Pi _i$, $i=1,2$, is
constant along $\partial M$ in each one of the two support hyperplanes, can be written as the fact that $\partial u / \partial \eta$ is
constant along the boundary $\partial M$. We can replace this boundary condition by a more general one and still obtain the same conclusion.

\begin{theorem}
\label{tmediabordo}
{\it
Let $\Pi _1$ and $\Pi _2$ be two parallel hyperplanes in $\realnum$.
Consider $M$ an embedded compact connected $C^2$ hypersurface contained in the region of $\realnum$ between  $\Pi _1$ and $\Pi _2$, with
$\partial D _i = \partial M \cap \Pi _{i} \neq \emptyset$, $i=1,2$. Suppose that the mean curvature of $M$ depends
only on the distance $u$ to $\Pi _i$, $i=1,2$, and $\displaystyle{\partial u / \partial \eta =h_{i}(H_{0i})}$, on $\partial D_i$, where $h_i$ is a continuous differentiable nonincreasing function of the mean curvature $H_{0i}$ of the boundary $\partial D_i$. Then $M$ is rotationally
symmetric with respect to a perpendicular orthogonal to $\Pi _i$, $i=1,2$ .}
\et

{\bf Proof}  We will begin repeating the same steps as in the proof of the first theorem: Let $L$ be a line in
$\Pi_1$ or $\Pi_2$ and $\Pi$ be a hyperplane perpendicular to $L$ and as a consequence also perpendicular to $\Pi_1$ and $\Pi_2$.
     Exists a hyperplane parallel to $\Pi$, that is tangent to $M$ or touches $M$ on its boundary and $M$ is entirely on one side of
     $\Pi$. We denote such hyperplane by $\Pi (0)$.
%
     By using the pushing argument we find a
     first contact point  $P$.
     If $P$ is an interior point of $\overline{M} (t^*)$,  we can prove that this hyperplane is a symmetry plane for $M$,
     following the same steps as in section 3.

     In the case when $P$ is a boundary point of $\overline{M} (t^*)$, we need a little more care.     As before,
     we choose a coordinate system $(x_1,...,x_n,u)$ in $\realnum$, with the origin at $P$ such that the tangent space to $M$ at $P$ is $u=0$.
     We can also choose the direction of the $u$ axis to be directed to the interior of $\Omega$, the direction of $x_1$ to the interior of $M$,
     so that the tangent space to $\partial M$ at $P$ contained in $\Pi _1$ or $\Pi_2$ is given by $u=x_1=0$.

     We have $M$ a $C^2$ hypersurface with boundary, so in a neighborhood of $P$,   $M$ and $\overline{M} (t^*)$ can be represented
     nonparametrically by $C^2$ functions $u(x)$ and $\overline{u}(x)$, defined on a domain $B \subset \realn$. Both functions satisfy:
     \begin{itemize}
      \item At the origin $u = \overline{u} $;
      \item On $\overline{B}$, $u(x)$ and $\overline{u} (x)$ satisfy equation of the mean curvature $H(x,u)$ and $u(x) \leq \overline{u} (x)$ because $\overline{M} (t)$ still lies in the interior of $\Omega$, for $t< t^*$;
      \item At the origin $\displaystyle{\frac{\partial u}{\partial \eta} \leq
         \frac{\partial \overline{u}}{\partial \eta}}$;
      \end{itemize}

      Let $H_{0i}$ be the mean curvature of $\partial M$ in $P$ and $\overline{H_{0}}$ be the mean curvature of $\partial \overline{M} (t^*)$ in $P$. Again, since
      $\overline{M} (t)$ still lies in the interior of $\Omega$, for $t< t^*$, $H_{0i} \leq \overline{H_0}$ in $P$. Using the hypothesis that $h_i$  is
      nonincreasing, we have at the origin
      \begin{equation*}
      \frac{\partial u}{\partial \eta} = h_{i}(H_{0i}) \geq h(\overline{H_0}) =
         \frac{\partial \overline{u}}{\partial \eta}\;,
      \end{equation*}
         for $i=1$ or $i=2$.

     Then, at the origin, $\displaystyle{\frac{\partial u}{\partial \eta} =
         \frac{\partial \overline{u}}{\partial \eta}}$ and we can conclude that $u(x)\equiv \overline{u} (x)$.

     As usual, we change the direction of $L$ and get to the fact that $M$ is rotationally
     symmetric with respect to a perpendicular orthogonal to $\Pi _i$, $i=1,2$.

     \hfill $\square$

\medskip

\begin{theorem}
\label{tbororadial}
{\it
Let $\Pi _1$ and $\Pi _2$ be two parallel hyperplanes in $\realnum$.
 Suppose that the origin $O$
 of $\realnum$ is inside the body $\Omega$. Consider $M$ an embedded compact connected $C^2$ hypersurface contained in the region of $\realnum$ between  $\Pi _1$ and $\Pi _2$, with $\partial M \cap \Pi _{i} \neq \emptyset$, $i=1,2$. Suppose that the mean curvature of $M$ depends
only on the distance $u$ to $\Pi _i$, $i=1,2$, and $\partial u / \partial \eta =-cr$, on the boundary of $M$, with $c>0$  a constant and $r$  the distance of a point to the origin. Then $M$ is rotationally
symmetric with respect to a perpendicular orthogonal to $\Pi _i$, $i=1,2$, containing the origin.}
\et

{\bf Proof}  Let $L$ be a line passing through the origin $O$ and contained in a parallel plane to $\Pi_1$ and $\Pi_2$.
Let $\Pi$ be a hyperplane perpendicular to $L$ and as a consequence also perpendicular to $\Pi_1$ and $\Pi_2$.
     Then there exists a hyperplane parallel to $\Pi$ that is tangent to $M$ or touches $M$ on its boundary and $M$ is entirely on one side of
     $\Pi$. We denote such hyperplane by $\Pi (0)$. Applying the pushing argument, we find a first touching point $P$.

     If $P$ is an interior point of $\overline{M} (t^*)$,  we can prove that this hyperplane is a symmetry plane for $M$,
     in a similar way as  in section 3. Since on the boundary of $M$ the values of $\partial u / \partial \eta =-cr$ have to coincide, the hyperplane
      $\Pi (t^*)$ has to be the one with contains the origin $O$.

     In the second case, if $P$ is a boundary point of $\overline{M} (t^*)$, we need again a little more care.

     By hypothesis, $\partial u / \partial \eta =-cr$.
     We choose a coordinate system $(x_1,...,x_n,u)$ in $\realnum$, with the origin at $P$ such that the tangent space to $M$ at $P$ is $u=0$. We can also choose the direction of the $u$ axis to be directed to the interior of $\Omega$, the direction of $x_1$ to the interior of $M$, so that the tangent space to $\partial M$ at $P$ contained in $\Pi _1$ or $\Pi_2$ is given by $u=x_1=0$.

     We have $M$ a $C^2$ hypersurface with boundary, so in a neighborhood of $P$,   $M$ and $\overline{M} (t^*)$ can be represented
     nonparametrically by $C^2$ functions $u(x)$ and $\overline{u}(x)$, defined on a domain $B \subset \realn$. Both functions satisfy:
     \begin{itemize}
      \item At the origin $u = \overline{u} $;
      \item On $\overline{B}$, $u(x)$ and $\overline{u} (x)$ satisfy equation of the mean curvature $H(x,u)$ and $u(x) \leq \overline{u} (x)$ because $\overline{M} (t)$ still lies in the interior of $\Omega$, for $t< t^*$.
      \item  At the origin, $\displaystyle{\frac{\partial u}{\partial \eta}=\frac{\partial \overline{u}}{\partial \eta}}$, because $M$ and $\overline{M} (t^*)$ have the same tangent plane at $P$, since it´s only possible to find one point of contact on the boundary.
     \end{itemize}

     Then, the plane of symmetry has to be the one with contains the origin and $u(x)\equiv \overline{u} (x)$.

     As usual, we change the direction of $L$ and get to the fact that $M$ is rotationally
     symmetric with respect to a perpendicular orthogonal to $\Pi _i$, $i=1,2$, containing the origin. \hfill$\square$

\section{Final Conclusions}
In \cite{athana} and \cite{vogel1}, M. Athanassenas and Vogel treated
the following free boundary problem: Find the orientable, compact
surface, embedded in $\R^3$, of least area and enclosing a fixed
volume, contained between two parallel planes and subject to the
condition that the boundary is constrained to lie on these two
parallel planes, studying the stability behavior of solutions.

A work by the present author is under preparation, investigating the stability
when we include the influence of gravity, in the same dimension as the previous
articles with constant mean curvature.

\medskip

\noindent INSTITUTO DE MATEM\'ATICA, UNIVERSIDADE FEDERAL DO RIO DE JANEIRO, ILHA DO FUND\~{A}O, 21941-909, RIO DE JANEIRO, BRAZIL

\noindent {\it E-mail address:}  {\bf monica@im.ufrj.br}

\end{document}